\documentclass[12pt]{article}
\usepackage{amssymb}
\usepackage{mathrsfs}
\usepackage{amsmath,amssymb,theorem}
\usepackage{graphicx}
\usepackage{subfigure}
\usepackage{psfrag}
\usepackage{color}
\newtheorem{thm}{Theorem}[section]
\newtheorem{conj}[thm]{Conjecture}
\newtheorem{cor}[thm]{Corollary}
\newtheorem{lem}[thm]{Lemma}
\theorembodyfont{\rmfamily}
\newtheorem{defn}[thm]{Definition}
\def\pf{\bigskip\noindent {\bf Proof.}~~}

\def\less{\backslash}

\def\pf{\bigskip\noindent {\emph{Proof.}}~~}

\def\qed{ \hfill $\square$}

\title{The Tur\'{a}n number of $2P_{7}$}

\author{\small {Yongxin Lan$^1$, Zhongmei Qin$^2$\footnote{The corresponding author.}\ ~and Yongtang Shi$^1$}\\
{\small  $^1$Center for Combinatorics and LPMC}\\
{\small Nankai University, Tianjin 300071, P.R. China}\\
{\small $^2$College of Science}\\
{\small Chang'an University, Xi'an, Shaanxi 710064, P.R. China}\\
{\small Emails: lan@mail.nankai.edu.cn, qinzhongmei90@163.com, shi@nankai.edu.cn}\\
}
\date{}
\begin{document}\maketitle
\begin{abstract}
The Tur\'an number of a graph $H$, denoted by $ex(n,H)$, is the maximum number of edges in any graph on $n$ vertices which does not contain $H$ as a subgraph. Let $P_{k}$ denote the path on $k$ vertices and let $mP_{k}$ denote $m$ disjoint copies of $P_{k}$. Bushaw and Kettle [Tur\'{a}n numbers of multiple paths and equibipartite forests, Combin. Probab. Comput. 20(2011) 837--853] determined the exact value of $ex(n,kP_\ell)$ for large values of $n$. Yuan and Zhang [The Tur\'{a}n number of disjoint copies of paths, Discrete Math. 340(2)(2017) 132--139] completely determined the value of $ex(n,kP_3)$ for all $n$, and also determined $ex(n,F_m)$, where $F_m$ is the disjoint union of $m$ paths containing at most one odd path. They also determined the exact value of $ex(n,P_3\cup P_{2\ell+1})$ for $n\geq 2\ell+4$. Recently, Bielak and Kieliszek [The Tur\'{a}n number of the graph $2P_5$, Discuss. Math. Graph Theory 36(2016) 683--694], Yuan and Zhang [Tur\'{a}n numbers for disjoint paths, arXiv: 1611.00981v1] independently determined the exact value of $ex(n,2P_5)$. In this paper, we show that $ex(n,2P_{7})=\max\{[n,14,7],5n-14\}$ for all $n \ge 14$, where $[n,14,7]=(5n+91+r(r-6))/2$, $n-13\equiv r\,(\text{mod }6)$ and $0\leq r< 6$. \\[2mm]
\end{abstract}

\section{Introduction}
Throughout this paper, we only consider simple graphs. For a graph $G$ we use $V(G)$, $|G|$, $E(G)$, $e(G)$ to denote the vertex set, number of vertices, edge set, number of edges, respectively. For $S_1,S_2\subseteq V(G)$ and $S_1\cap S_2=\emptyset$, denote by $e(S_1,S_2)$ the
number of edges between $S_1$ and $S_2$. Let $G$ and $H$ be two disjoint graphs. By $G\cup H$ denote the disjoint
union of graphs $G$ and $H$ and by $kG$ denote the $k$ disjoint
copies of $G$. Denote by $G+H$ the graph obtained from $G\cup
H$ by joining all vertices of $G$ to all vertices of
$H$. Let $\overline{G}$ be the complement of the graph $G$.
Denote by $P_n$, $C_n$ and
$K_{n}$ the path, cycle and complete graph on $n$ vertices, respectively. For $S\subseteq V(G)$, let $G[S]$ denote the subgraph of $G$ induced by
$S$ and let $|S|$ denote the cardinality of $S$. For a graph $G$ and its subgraph $H$, by $G-H$ we mean a graph obtained from $G$ by deleting all vertices of $H$ with all incident edges. If $H$ consists of a single vertex $x$, then we simple write $G-x$. For $v\in V(G)$, let
$N_{G}(v)$ denote the set of vertices in $G$ which are adjacent to $v$. We define $d_{G}(v)=|N_G(v)|$.

A graph is {\it $H$-free} if it does not contain $H$ as a subgraph. The $Tur\acute{a}n$ $number$ of a graph $H$, denoted by $ex(n,H)$, is the
maximum number of edges in any $H$-free graph on $n$ vertices, i.e., $$ex(n,H)=\max\{e(G): H\nsubseteqq G \text{ and } |G|=n\}.$$
Let $EX(n,H)$ denote the family of all $H$-free graphs on $n$ vertices with $ex(n,H)$ edges. A graph in $EX(n,H)$ is called an {\it extremal graph} for $H$.
Moreover, we denote by $ex_{con}(n,H)$ the maximum number of edges
in any connected $H$-free graph on $n$ vertices. 
The problem of determining Tur\'{a}n number for
assorted graphs traces its history back to 1907, when Mantel (see,
$e.g.$, \cite{BollobMGT}) proved  $ex(n,C_3)=\lfloor n^{2}/4\rfloor$.
In $1941$, Tur\'an \cite{Tur41,Tur54} proved that the
extremal graph for $K_{r}$ is the complete
$(r-1)$-partite graph, which is as balanced as
possible (any two part sizes differ at most 1). The balanced complete
$(r-1)$-partite graph on $n$ vertices is called as the Tur\'{a}n graph denoted by
$T_{r-1}(n)$. For sparse graphs, Erd\H{o}s and Gallai \cite{ErdGal} in
$1959$ proved the following well known result.
\begin{thm}[\cite{ErdGal}]
Let $G$ be a $P_k$-free graph on $n$ vertices and $n\geq k\geq 2$. Then $e(G)\leq
(k-2)n/2$ with equality if and only if $n=(k-1)t$ and
$G=tK_{k-1}$.
\end{thm}

For convenience, we first introduce the following symbols.
\begin{defn}
Let $n\geq m\geq\ell\geq 3$ be given three positive integers. If
$n$ can be written as $n=(m-1)+t(\ell-1)+r$, where $t\geq 0$ and $0\leq
r< \ell-1$, then we denote
$$[n,m,\ell]= {m-1 \choose 2}+t{\ell-1 \choose 2}+{r \choose 2}.$$
Moreover, if $n\leq m-1$, then we denote
$$[n,m,\ell]= {n \choose 2}.$$
\end{defn}
\begin{defn}
Let $s=\sum_{i=1}^{m}\lfloor k_i/2\rfloor$ and $k_{i}$
be positive integers. If $n\ge s$, then we denote
$$[n,s]={s-1 \choose 2}+(s-1)(n-s+1).$$
\end{defn}

Later, for all integers $n$ and $k$, Faudree and Schelp \cite{FaudSche} characterized all extremal graphs for $P_k$.
\begin{thm}[\cite{FaudSche}]\label{ALLPATH}
Let $G$ be a graph on $n=t(k-1)+r$ $(0\leq t$ and $0\leq r<k-1)$ vertices. If $G$ is $P_k$-free, then $e(G)\leq [n,k,k]$. Moreover, the equality  holds if
and only if

$\bullet$ $G=(tK_{k-1})\cup
K_{r}$ or
\medskip{}

$\bullet$ $G=((t-s-1)K_{k-1})\cup(K_{(k-2)/2}+\overline{K_{k/2+s(k-1)+r}})$, where $k$ is even, $t>0$, $r=k/2$ or $(k-2)/2$ and $0\leq s<t$.
\end{thm}
\begin{cor}
For a positive integer $n\equiv r\,(mod~k)$, $ex(n,P_{k+1})=(n(k-1)+r(r-k))/2.$
\end{cor}

We see that $ex(n,P_{k})$ has been determined for all integers $n \ge k$ and
all extremal graphs has also been characterized. For connected
graphs, Kopylov \cite{Kopy} and  Balister, Gy\H{o}ri, Lehel, and Schelp \cite{BalisGy} determined
$ex_{con}(n,P_{k})$ and characterized all extremal graphs for all integers $n \ge k$. Recently, Lan, Shi and Song \cite{LSS17} studied the Tur\'an number of paths in planar graphs.
\begin{thm}[\cite{BalisGy,Kopy}] \label{CONPATH}
Let $G$ be a connected $P_k$-free graph on $n$ vertices and $n\geq k\geq 4$. Then \\
$$e(G)\leq \max\left\{{k-2 \choose 2}+(n-k+2),[n,\lfloor k/2\rfloor]+c\right\},$$
where $k\equiv c~(mod~2)$. Further, the equality holds if and only if $G=(K_{k-3}\cup \overline{K_{n-k+2}})+K_{1}$
or
$G=(K_{1+c}\cup
\overline{K_{n-\lfloor(k+1)/2\rfloor}})+K_{\lfloor k/2\rfloor-1}.$
\end{thm}

In 1962, Erd\H{o}s \cite{Erdos} first studied the Tur\'{a}n
number of $kK_{3}$. And later, Moon \cite{Moon} and Simonovits \cite{Simonovits} studied the case of $kK_{r}$. In 2011, Bushaw
and Kettle \cite{Bushaw} determined $ex(n,kP_\ell)$ for $n$
sufficiently large.

\begin{thm}[\cite{Bushaw}]
For integers $k\geq2$, $\ell\geq4$ and $n\geq 2\ell+2k\ell(\lceil\ell/2\rceil+1){\ell
\choose \lfloor\ell/2\rfloor}$, $$ex(n,k
P_{\ell})=\left[n,k\left\lfloor\frac\ell 2\right\rfloor\right]+c,$$
where $\ell\equiv c~(mod~2)$.
\end{thm}

Furthermore, their proof shows that their construction is optimal for
$n=\Omega(k\ell^{3/2}2^{\ell})$. Moreover, Bushaw and Kettle
conjectured that their construction is optimal for $n=\Omega(k\ell)$. Recently, Lidick\'{y} et al. \cite{LidickLiuPalmer} extended Bushaw and Kettle's result and
determined $ex(n,F_{m})$ for $n$ sufficiently large, where $F_{m}=\bigcup_{i=1}^{m}P_{k_{i}}$
and $k_{1}\geq k_{2}\geq\ldots\geq k_{m}$.
\begin{thm}[\cite{LidickLiuPalmer}]
Let $F_{m}=\bigcup_{i=1}^{m}P_{k_{i}}$ and $k_{1}\geq k_{2}\geq\ldots\geq k_{m}$. If at least
one $k_{i}$ is not 3, then for $n$ sufficiently large,
$$ex(n,F_{m})=\left[n,\sum_{i=1}^{m}\left\lfloor\frac{k_{i}}{2}\right\rfloor\right]+c,$$
where $c=1$ if all $k_{i}$ are odd, and $c=0$ otherwise. Moreover,
the extremal graph is unique.
\end{thm}

However, they did not consider $ex(n,F_{m})$ for
smaller $n$. Recently,
Yuan and Zhang \cite{ZhangXiaodong,kP3} completely determined the value of $ex(n,kP_3)$ and characterized all the extremal graphs for all $n$. Furthermore, they proved the following result in which $F_m$ contains at most one odd path and proposed Conjecture \ref{Conjecture}.
\begin{thm}[\cite{ZhangXiaodong}]
Let $k_{1}\geq k_{2}\geq \ldots \geq k_{m}\geq3$, $n\geq
\sum_{i=1}^{m}k_{i}$ and $F_{m}=\bigcup_{i=1}^{m}P_{k_{i}}$. If there is at most one odd in
$\{k_{1},k_{2},\ldots,k_{m}\}$, then
$$ex(n,F_{m})=\max\left\{[n,k_{1},k_{1}],[n,k_{1}+k_{2},k_{2}],\ldots,\left[n,\sum_{i=1}^{m}k_{i},
k_{m}\right],\left[n,\sum_{i=1}^{m}\left\lfloor\frac{k_{i}}{2}\right\rfloor\right]\right\}.$$
Moreover, if $k_{1},k_{2},\ldots,k_{m}$ are even, then the extremal
graphs are characterized.
\end{thm}
\begin{conj}[\cite{ZhangXiaodong}] \label{Conjecture}
Let $k_{1}\geq k_{2}\geq\ldots\geq k_{m}\geq3$, $k_{1}>3$ and
$F_{m}=\bigcup_{i=1}^{m}P_{k_{i}}$. Then
$$ex(n,F_{m})=\max\left\{[n,k_{1},k_{1}],[n,k_{1}+k_{2},k_{2}],\ldots,\left[n,\sum_{i=1}^{m}k_{i},
k_{m}\right],\left[n,\sum_{i=1}^{m}\left\lfloor\frac{k_{i}}{2}\right\rfloor\right]+c\right\},$$
where $c=1$ if all of $k_{1},k_{2},\ldots,k_{m}$ are odd, and $c=0$
otherwise. Moreover, the extremal graphs are $$EX(n,P_{k_{1}}),\ldots,K_{\sum_{i=1}^{m}k_{i}-1}\cup H ~\text{for}~H\in EX(n-\begin{matrix}\sum_{i=1}^{m}k_{i}+1\end{matrix},P_{k_{m}}), ~and$$
$$K_{\sum_{i=1}^{m}\lfloor k_{i}/2\rfloor-1}+(K_{1+c}\cup
\overline{K_{n-\sum_{i=1}^{m}\lfloor k_{i}/2\rfloor-c}}).$$
\end{conj}

When there are at least two odd integers in $\{k_{1},k_{2},\ldots,k_{m}\}$, Yuan and Zhang also determined $ex(n,P_3\cup P_{2\ell+1})$ for $n\geq 2\ell+4$ and characterized all extremal graphs. Bielak and Kieliszek \cite{2P5} and Yuan and Zhang \cite{ZhangXiaodong} independently determined $ex(n,2P_5)$ and characterized all extremal graphs. In this paper, we prove the following result, which partially confirms Conjecture \ref{Conjecture}.
\begin{thm}\label{them}
For $n\geq 14$, $$ex(n,2P_{7})=\max\{[n,14,7],5n-14\}.$$
Moreover, the extremal graphs are $K_{13}\cup H$ for $H\in EX(n-13,P_{7})$ when $n\le 21$ and
$K_{5}+(K_{2}\cup \overline{K_{n-7}})$ when $n\ge 22$.
\end{thm}

\section{Proof of Theorem \ref{them}}
We first present some useful lemmas. In the following, we say that $u$ {\it hits} $v$ or $v$ {\it hits} $u$ if two vertices $u$ and $v$ are adjacent. Otherwise, we say that $u$ {\it misses} $v$ or $v$ {\it misses} $u$ if $u$ and $v$
are not adjacent. We say a vertex set $A$ hits (misses) a vertex set $B$, it means that each vertex of $A$ is adjacent (non-adjacent) to each vertex of $B$.
\begin{lem}[Observation 2 of \cite{ZhangXiaodong}]\label{NKK}
Let $k_{1}\geq k_{2}\geq 3$ be two positive integers. If $n_{1}\geq
k_{1}$, then $[n_{1},k_{1}+k_{2},k_{2}]+[n_{2},k_{2},k_{2}]\leq
[n_{1}+n_{2},k_{1}+k_{2},k_{2}].$
\end{lem}
\begin{lem}[Observation 5 of \cite{ZhangXiaodong}]\label{Nk}
Let $k_{1}\geq k_{2}\geq 3$ be two positive integers. If $n_{1}\geq
k_{1}+k_{2}$, then
$[n_{1},\lfloor k_1/2\rfloor+\lfloor k_2/2\rfloor]+[n_{2},k_{2},k_{2}]<[n_{1}+n_{2},\lfloor k_1/2\rfloor+\lfloor k_2/2\rfloor].$
\end{lem}
\begin{lem}\label{N14}
Let $G$ be a connected $2P_7$-free graph on $n\geq 14$ vertices. Then $$e(G)\le\max\{[n,14,7],5n-14\}.$$
with equality only when $n\ge 22$ and $G=K_{5}+(\overline{K_{n-7}}\cup K_{2})$.
\end{lem}
\pf
Let $G\neq K_{5}+(\overline{K_{n-7}}\cup K_{2})$ be any connected $2P_7$-free
graph on $n$ vertices with $e(G)\geq \max\{[n,14,7],5n-14\}$ edges. Note that $\max\{[n,14,7],5n-14\}=[n,14,7]$ when $n\le 21$ and $\max\{[n,14,7],5n-14\}=5n-14$ when $n\ge 22$. Since $\max\{[n,14,7],5n-14\}\geq
ex_{con}(n,P_{13})$, by Theorem \ref{CONPATH}, $G$ contains
$P_{13}$ as a subgraph.
Let $P_{13}=x_{1}x_{2}\ldots x_{13}$ be a subgraph of $G$. Then\medskip

\noindent($\ast$) each vertex of $G-P_{13}$
cannot hit two adjacent vertices in $P_{13}$.\medskip

Notice that each vertex in $G-P_{13}$ misses $\{x_{1},x_{6},x_{8},x_{13}\}$
and can not hit both $x_{p}$ and $x_{p+8}$ for $p\in \{2,3,4\}$. Moreover, if $y$ is an isolated vertex in $G-P_{13}$, then by ($\ast$), $|N_G(y)\cap V(P_{13})|\leq5$; if $y$ is not an
isolated vertex in $G-P_{13}$, then $N_G(y)\cap V(P_{13})\subseteq\{x_{3},x_{4},x_{7},x_{10},x_{11}\}$ and so $|N_G(y)\cap V(P_{13})|\leq3$ by ($\ast$); if $P_{k}=y_{1}y_{2} \ldots y_{k}\subseteq
G-P_{13}$ and $k\geq 3$ such that $y_{1}$ hits $P_{13}$, then
$y_{1}$ can only hit $x_{7}$. Now we will prove the following
useful Facts.\\

\noindent\textbf{Fact 1.} $e(G[V(P_{13})])\le 74$.\medskip

Since $G$ is connected and $n\ge 14$, at least one vertex of $V(G)\less V(P_{13})$ hits $P_{13}$, say $x_i$. Then either $i\ge 6$ or $i\le 8$. Without loss of generality, we may assume that $i\ge 6$. For $1\le j\le i-2$, if both $x_{13}x_{j}\in E(G)$ and $x_{i+1}x_{j+1}\in E(G)$, then $G$ contains $2P_7$ as a subgraph, a contradiction. Thus $e(G[V(P_{13})])\le 74$.\medskip

\noindent\textbf{Fact 2.} If there exists a $P_{3}=y_{1}y_{2}y_{3}\subseteq G-P_{13}$ such that
$y_{1}$ hits $P_{13}$, then we have $e(G[V(P_{13})])\leq 57.$\medskip

Clearly, $y_{1}$ must hit $x_{7}$ and so $G$ contains a copy of $P_7$ with vertices $x_4,x_5,x_6,x_7,y_1,y_2,y_3$. Therefore, $\{x_1,x_2,x_3,x_5,x_6\}$ misses
$\{x_{11},x_{12},x_{13}\}$. Symmetrically, $\{x_8,x_9,x_{11},x_{12},x_{13}\}$ misses $\{x_1,x_2,x_3\}$.
So $e(G[V(P_{13})])\leq 78-(2\cdot15-9)=57$.\qed\\

\noindent\textbf{Fact 3.}
If there exists a non-isolated vertex in $G-P_{13}$, that hits one vertex of $P_{13}$, then we have $e(G[V(P_{13})])\leq 68.$\medskip

Let $y$ be a non-isolated vertex in $G-P_{13}$, that hits one vertex, say $x_i$ of $P_{13}$. Recall that $x_i\in \{x_3,x_4,x_7,x_{10},x_{11}\}$. If $x_i\in \{x_3,x_4\}$, then $\{x_{1},x_{2},\dots,x_{i-1}\}$ misses
$\{x_{i+1},x_{i+2},x_{9},x_{12},x_{13}\}$ and so
$e(G[V(P_{13})])\leq 68$. Symmetrically, if $x_i\in \{x_{10},x_{11}\}$, then
$e(G[V(P_{13})])\leq 68$.
Now assume that $x_i=x_7$. Then $\{x_{1},x_{2},x_{i-1},x_{i-2}\}$ misses
$\{x_{12},x_{13}\}$ and symmetrically $\{x_{i+1},x_{i+2},x_{12},x_{13}\}$ misses
$\{x_{1},x_{2}\}$. So $e(G[V(P_{13})])\leq 78-(2\cdot 8-4)=66$.\qed\\

\noindent\textbf{Fact 4.}
If there exists a non-isolated vertex in $G-P_{13}$, that hits two vertices of $P_{13}$, then we have $e(G[V(P_{13})])\leq 59$.\medskip

Let $y$ be a non-isolated vertex in $G-P_{13}$, that hits two vertices, say $x_i$ and $x_j$ ($i<j$), of $P_{13}$. Recall that $\{x_i,x_j\}\subseteq \{x_3,x_4,x_7,x_{10},x_{11}\}$ and $\{x_i,x_j\}\ne \{x_3,x_{11}\}$. If $x_i=x_3$, then by ($\ast$), $x_j\in\{x_7,x_{10}\}$. Thus $\{x_{1},x_{2}\}$ misses
$\{x_{4},x_{5},x_{6},x_{8},x_{9},x_{11},x_{12},x_{13}\}$ and
$\{x_{j-2},x_{j-1}\}$ misses $\{x_{12},x_{13}\}$. So $e(G[V(P_{13})])\leq 58$. Symmetrically,
if $x_j=x_{11}$, then by ($\ast$), $x_i\in\{x_4,x_7\}$ and so $e(G[V(P_{13})])\leq 58$. Now we can assume that $x_i\neq x_3$ and $x_j\neq x_{11}$.
If $x_i=x_4$, then $x_j\in\{x_7,x_{10}\}$. Thus $\{x_{1},x_{2},x_3\}$ misses
$\{x_{5},x_{6},x_{9},x_{12},x_{13}\}$ and $\{x_{j-2},x_{j-1}\}$ misses $\{x_{12},x_{13}\}$. So $e(G[V(P_{13})])\leq 59$. Symmetrically,
if $x_j=x_{10}$, then $x_i\in\{x_4,x_7\}$ and so $e(G[V(P_{13})])\leq 59$.\qed\\

\noindent\textbf{Fact 5.}
If there exists an isolated vertex in $G-P_{13}$, that hits five vertices of
$P_{13}$, then $e(G[V(P_{13})])\leq 50$.\medskip

Let $y$ be an isolated vertex in $G-P_{13}$, that hits exactly five
vertices, say $x_i, x_j, x_k, x_\ell, x_m$, $i<j<k<\ell<m$ of $P_{13}$. Recall that $\{x_i, x_j, x_k, x_\ell, x_m\}\subseteq V(P_{13})\setminus \{x_1,x_6,x_8,x_{13}\}$ and $y$ cannot hit both $x_{p}$ and $x_{p+8}$ for $p\in \{2,3,4\}$. Since $y$ cannot hit two adjacent vertices in $P_{13}$, we have $x_k=x_7$, $\{x_i,x_j\}\subseteq \{x_{2},x_{3},x_{4},x_{5}\}$ and $\{x_\ell,x_m\}\subseteq \{x_{9},x_{10},x_{11},x_{12}\}$. Let $A=\{x_{i-1},x_{j-1},x_{k-1},x_{\ell-1},x_{m-1},x_{13}\}$ and $B=\{x_1,x_{i+1},x_{j+1},x_{k+1},x_{\ell+1},x_{m+1}\}$. Then,
$A$ and $B$ are independent sets and $|A\cap B|=4$. Since $\{x_3,x_{11}\}\nsubseteqq N_G(y)$, we have either $i=2$ or $m=12$. 
If $i=2$ and $m=12$, then $N_G(y)=\{x_2,x_5,x_7,x_9,x_{12}\}$, which implies that $x_5$ misses $\{x_{10},x_{11}\}$. And symmetrically $x_9$ misses $\{x_3,x_4\}$. If $i=2$ and $m\neq12$, then $\ell=9$ and $m=11$, which implies that $x_m$ misses $\{x_3,x_6\}$ and $x_\ell$ misses $\{x_q,x_{q+1}\}\subseteq \{x_1,\dots,x_7\}\setminus N_G(y)$. If $i\ne 2$ and $m=12$, then $i=3$ and $j=5$, which implies that $x_i$ misses $\{x_8,x_{11}\}$ and $x_j$ misses $\{x_q,x_{q+1}\}\subseteq \{x_7,\dots,x_{13}\}\setminus N_G(y)$. For each of the above cases, we have $e(G[V(P_{13})])\leq 78-({|A| \choose 2}+{|B| \choose 2}-{|A\cap B| \choose 2})-4=50$.\qed\\

\noindent\textbf{Fact 6.}
If there exists an isolated vertex in $G-P_{13}$, that hits four vertices of
$P_{13}$, then $e(G[V(P_{13})])\leq 59$.\medskip

Let $y$ be an isolated vertex in $G-P_{13}$, that hits exactly four
vertices, say $x_i, x_j, x_k, x_\ell$, $i<j<k<\ell$ of $P_{13}$. Recall that $\{x_i, x_j, x_k, x_\ell\}\subseteq V(P_{13})\setminus \{x_1,x_6,x_8,x_{13}\}$ and $y$ cannot hit both $x_{p}$ and $x_{p+8}$ for $p\in \{2,3,4\}$. Let $A=\{x_{i-1},x_{j-1},x_{k-1},x_{\ell-1},x_{13}\}$ and $B=\{x_1,x_{i+1},x_{j+1},x_{k+1},x_{\ell+1}\}$. Then $A$ and $B$ are independent sets and $|A\cap B|\le 3$. If $|A\cap B|\le 2$, then $e(G[V(P_{13})])\leq 78-({|A| \choose 2}+{|B| \choose 2}-1)=59$. Now we assume that $|A\cap B|=3$. If $i=2$ and $\ell=12$, then $7\in\{j,k\}$ which implies that $x_3$ misses $x_{11}$ and $x_p$ misses $x_{p+9}$ for $p\in \{1,4\}$. If $i=2$, $\ell\ne 12$ and $7\in\{j,k\}$, then $x_{11}$ misses $\{x_3,x_6\}$. If $i=2$, $\ell\ne 12$ and $7\notin\{j,k\}$, then $N_G(y)=\{x_2,x_4,x_9,x_{11}\}$ which implies $x_{11}$ misses $\{x_5,x_8\}$. If $\ell=12$ and $i\ne 2$, then it is similar as the case of $i=2$ and $\ell\ne 12$. If $i\ne 2$ and $\ell\ne 12$, then $N_G(y)=\{x_3,x_5,x_7,x_9\}$ which implies $x_{11}$ misses $\{x_1,x_4\}$. For each of the above cases, $e(G[V(P_{13})])\leq 78-({|A| \choose 2}+{|B| \choose 2}-{|A\cap B| \choose 2})-2=59$.\qed\\

\noindent\textbf{Fact 7.}
If there exists an isolated vertex in $G-P_{13}$, that hits three vertices of
$P_{13}$, then $e(G[V(P_{13})])\leq 67$.\medskip

Let $y$ be an isolated vertex in $G-P_{13}$, that hits exactly three
vertices, says $x_i, x_j, x_k$, $i<j<k$ of $P_{13}$. Recall that $\{x_i, x_j, x_k\}\subseteq V(P_{13})\setminus \{x_1,x_6,x_8,x_{13}\}$ and $y$ can not hit both $x_{p}$ and $x_{p+8}$ for $p\in \{2,3,4\}$. Let $A=\{x_{i-1},x_{j-1},x_{k-1},x_{13}\}$ and $B=\{x_1,x_{i+1},x_{j+1},x_{k+1}\}$. Then both $A$ and $B$ are independent sets and $|A\cap B|\le 2$. Hence, $e(G[V(P_{13})])\leq 78-({|A| \choose 2}+{|B| \choose 2}-{|A\cap B| \choose 2})\le 78-(6+6-1)=67$.\qed
\medskip{}

Let $P_{k}=y_{1}y_{2}\ldots y_{k}$, where $k\leq 6$, be the longest path in
$G-P_{13}$ such that $y_{1}$ hits $P_{13}.$ Let $H_1,H_2,\ldots,H_t$ be connected components of order at least 2 of $G-P_{13}$ and let $H$ be a subgraph of $G$ which consists of all isolated vertices of $G-P_{13}$. Note that $\sum_{i=1}^t|H_i|+|H|=n-13$. Let $m(H_i)$ be the number of edges incident with the vertices of $H_i$ and let $H_{1}$ be a component of $G-P_{13}$ which contains $P_{k}$ as a subgraph. We first show the following claim.
\medskip{}

\textbf{Claim:} For $1\leq i\leq t$, $m(H_i)\leq 4|H_i|$.

\begin{pf}
We use induction on $|H_i|$. Recall that each vertex of $H_i$ can hit at most three vertices of $P_{13}$. For $|H_i|=2$, $m(H_i)=e(G[V(H_i)])+e(V(H_i),V(P_{13}))\leq 7\leq4|H_i|$. If $H_i$ has a pendant vertex $x$, then $d_G(x)\leq 4$. By induction hypothesis, we have $m(H_i)=m(H_i-x)+d_G(x)\leq4(|H_i|-1)+4\leq4|H_i|$. Next if $H_i$ has no pendant vertex, then each vertex of $H_i$ must be an endpoint of a path of length at least two. This implies that each vertex of $H_i$ can only hit $x_{7}$ of $P_{13}$. Thus, $m(H_i)=e(G[V(H_i)])+e(V(H_i),V(P_{13}))\leq ex_{con}(|H_i|,P_7)+|H_i|\leq\frac{7}{2}|H_i|$ since $H_i$ is $P_7$-free.\qed
\end{pf}\\

Let $\Delta(H)=\max\{d_{G}(v)|v\in V(H)\}$. Recall that $\Delta(H)\le 5$. Now we would divide the proof into the following cases (in each case we assume, the previous cases do not hold).

\textbf{Case 1.} $\Delta(H)=5$. Then by Fact $5$ and the Claim,
$$e(G)\le 50+5(n-13)=5n-15<\max\{[n,14,7],5n-14\},$$ a contradiction.

\textbf{Case 2.} $\Delta(H)=4$ or $k\ge 3$ or there exists a non-isolated vertex in $G-P_{13}$ that hits two vertices of $P_{13}$ ($k=2$). Then by Facts $6$, $2$ and $4$ and the Claim,
$$e(G)\leq
59+4(n-13)=4n+7<\max\{[n,14,7],5n-14\},$$ a contradiction.

\textbf{Case 3.} $\Delta(H)=3$ ($k=2$) or there exists a non-isolated vertex in $G-P_{13}$ that hits one vertex of $P_{13}$ ($k=2$). For $k=2$, each component of $G-P_{13}$ is a star (with
at least three vertices), or an edge, or an isolated vertex. For $1\leq i\leq t$, $e(G[V(H_i)])\leq|H_i|-1$. $m_0\leq\sum_{i=1}^t(2|H_i|-1)+3|H|=3(n-15)+6-\sum_{i=1}^t|H_i|-t\leq3(n-13)$. Then by Facts $7$ and $3$, we have
$$e(G)\leq 68+3(n-13)=3n+29<\max\{[n,14,7],5n-14\},$$ a
contradiction.

\textbf{Case 4.} $\Delta(H)\le 2$ and $k=1$. Then by Fact $1$,
$$e(G)\leq 74+2(n-13)=2n+48<\max\{[n,14,7],5n-14\},$$ a
contradiction.

The proof is thus completed.\qed\\

\noindent{\it Proof of Theorem \ref{them}.}\
Let $G$ be any $2P_{7}$-free graph on $n$ vertices
with $e(G)\geq \max\{[n,14,7],\\5n-14\}$. If $G$ is connected, then by Lemma \ref{N14}, $e(G)\le \max\{[n,14,7],5n-14\}$ when $n\ge 22$ and $e(G)<\max\{[n,14,7],5n-14\}$ when $n\le 21$. Thus when $G$ is connected, $e(G)\le \max\{[n,14,7],5n-14\}$ with equality holds if and only if  $n\ge 22$ and $G=K_{5}+(\overline{K_{n-7}}\cup K_{2})$. Now we may assume that $G$ is disconnected. By Lemma \ref{ALLPATH}, $G$
contains $P_{7}$ as a subgraph. Let $C$ be a connected component with
$n_{1}\geq 7$ vertices which contains $P_{7}$ as a subgraph. Notice that $C$
is $2P_{7}$-free and $G-C$ is $P_{7}$-free.
If
$n_{1}\geq 22$, then by Lemma \ref{N14}, $e(C)\le 5n-14$ and by Lemmas \ref{ALLPATH} and \ref{Nk},
$$e(G)=e(C)+e(G-C)\leq 5n_{1}-14+[n-n_{1},7,7]<5n-14,$$
a contradiction. If $14\le n_1\le 21$, then by Lemma \ref{N14}, $e(C)<[n_1,14,7]$ and  by Lemmas \ref{ALLPATH} and \ref{NKK}, $$e(G)=e(C)+e(G-C)<[n_1,14,7]+[n-n_{1},7,7]\le [n,14,7],$$
a contradiction.
If $n_{1}\leq 13$, then $e(G)\leq {n_{1} \choose
2}+[n-n_{1},7,7]\leq [n,14,7]$ with equality holds if and only if $C=K_{13}$ and $G-C\in EX(n-13,P_{7})$. But then when $n\ge 22$, $e(G)\ge \max\{[n,14,7],5n-14\}=5n-14>[n,14,7]$, a contradiction. Thus when $G$ is disconnected, $e(G)\le \max\{[n,14,7],5n-14\}$ with equality holds if and only if  $n\le 21$, $G=K_{13}\cup H$ for $H\in EX(n-13,P_{7})$.

The proof is thus complete.\qed\\

\noindent{\bf Acknowledgements:} We wish to thank the two anonymous referees for their valuable
suggestions and comments. This work was supported by the National Science Foundation and the
Natural Science Foundation of Tianjin (No. 17JCQNJC00300).

\end{document}